\documentclass[11pt]{amsart}
\usepackage[all]{xy}

\usepackage{amsmath}
\usepackage{amsfonts}
\usepackage{amssymb}
\usepackage{mathrsfs}

\long\def\symbolfootnote[#1]#2{\begingroup%
\def\thefootnote{\fnsymbol{footnote}}\footnote[#1]{#2}\endgroup}

\DeclareFontEncoding{OT2}{}{} 
\newcommand{\textcyr}[1]{%
 {\fontencoding{OT2}\fontfamily{wncyr}\fontseries{m}\fontshape{n}
 \selectfont #1}}
\newcommand{\Sha}{{\mbox{\textcyr{Sh}}}}

\numberwithin{equation}{section}

\def\bz{{\mathbb Z}\,}
\def\bq{{\mathbb Q}}
\def\bg{{\mathbb G}}
\def\spec{{\rm{Spec}}\,}
\def\fb{\overline{F}}

\def\img{{\rm{Im}}\,}

\def\be{\kern -.1em}
\def\lbe{\kern -.025em}
\def\sho{\s H^{\le\circ}}

\def\s{\mathscr }
\def\ra{\rightarrow}
\def\e{\kern 0.08em}
\def\le{\kern 0.03em}
\def\ng{\kern -0.04em}

\def\g{\varGamma}
\def\gv{\varGamma_{\! v}}

\def\sao{{\s A}^{\lbe\circ}}
\def\sbo{{\s B}^{\le\circ}}
\def\krn{{\rm{Ker}}\,}
\def\cok{{\rm{Coker}}\,}

\def\sa{\s A}
\def\sov{\s O_{\be v}}
\def\gv{\varGamma_{\! v}}

\newtheorem{theorem}{Theorem 4.9\!\!}

\newtheorem{lemma}{Lemma}[section]
\newtheorem{teorema}[lemma]{Theorem}
\newtheorem{corollary}[lemma]{Corollary}
\newtheorem{proposition}[lemma]{Proposition}
\theoremstyle{definition}

\theoremstyle{remark}
\newtheorem{remark}[lemma]{Remark}
\newtheorem{remarks}[lemma]{Remarks}

\begin{document}

\title[N\'eron class groups of abelian varieties]{On N\'eron class groups of
abelian varieties}


\author{Cristian D. Gonz\'alez-Avil\'es}
\address{Departamento de Matem\'aticas, Universidad de La Serena, Chile}
\email{cgonzalez@userena.cl}


\thanks{The author is partially supported by Fondecyt grant
1080025}

\dedicatory{A la memoria de Esbriel Avil\'es Cruz.*}

\thanks{* See http://uls.cl/$\sim$cgonzalez/}

\maketitle

\begin{abstract} Let $F$ be a global field,
let $S_{\infty}$ be the set of archimedean primes of $F$ and let $S$
be any nonempty finite set of primes of $F$ containing $S_{\infty}$.
In this paper we study the N\'eron $S$-class group $C_{A,\e F,\e S}$
of an abelian variety $A$ defined over $F$. In the well-known
analogy that exists between the Birch and Swinnerton-Dyer conjecture
for $A$ over $F$ and the analytic class number formula for the field
$F$ (in the number field case), the finite group $C_{A,\e F,\e
S_{\infty}}$ (not the Tate-Shafarevich group of $A$) is a natural
analog of the ideal class group of $F$.
\end{abstract}

\section{Introduction}
Let $F$ be a global field, let $S$ be a nonempty finite set of
primes of $F$ containing all archimedean primes and let $\s O_{F,\e
S}$ be the ring of $S$-integers of $F$. Further, we write $U=\spec\s
O_{F,\e S}$. Note that $v\notin S$ if, and only if, $v$ corresponds
to a point of $U$.

A central problem in Number Theory is that of extending certain
known results for tori to abelian varieties over $F$. Specifically,
the analytic class number formula for a number field $F$ has long
been regarded as a template for the Birch and Swinnerton-Dyer
conjecture for an abelian variety $A$ over $F$. But the analogy
between a theorem for the trivial torus $T=\bg_{m,\e F}$ (as the
analytic class number formula certainly is) and a conjecture about
an arbitrary abelian variety is a distant one, and many researchers
have come to view (incorrectly, we believe) the Tate-Shafarevich
group $\!\!\Sha^{1}(A)$ of $A$ as a natural analog of the ideal
class group of $F$.

In this paper we introduce the {\it N\'eron $S$-class group}
$C_{A,\e F,\e S}$ of $A$ over $F$ and establish a duality theorem
for this group. The definition of $C_{A,\e F,\e S}$ is quite simple.
For each prime $v\in U $, let $F_{v}$ be the completion of $F$ at
$v$ and let $k(v)$ denote the corresponding residue field. There
exists a canonical {\it reduction map} $A(F)\ra\Phi_{v}(A)(k(v))$,
where $\Phi_{v}(A)$ is the group scheme of connected components of
the special fiber of the N\'eron model of $A_{F_{v}}$. Then
$$
C_{A,\e F,\e S}=\cok\!\!\left[\e A(F)\ra\bigoplus_{v\in U
}\Phi_{v}(A)(k(v))\right],
$$
where the points of $A(F)$ are mapped diagonally into
$\bigoplus_{\,v\in U }\Phi_{v}(A)(k(v))$ via the preceding reduction
maps. For reasons that are explained in Remark 3.3, we believe that
this group is the correct analog of the ideal $S$-class group of
$F$.

We will now state the main theorem of this paper. Let $B$ denote the
dual (i.e., Picard) variety of $A$. For each $v\in U $,
Grothendieck's pairing
$$
\Phi_{v}(A)\big(\e\overline{k(v)}\e\big)\times
\Phi_{v}(B)\big(\e\overline{k(v)}\e\big)\ra{\bq}/{\bz}
$$
induces a nondegenerate pairing of finite groups
$$
\Phi_{v}(A)(k(v))\times H^{\e 1}\lbe(k(v),\Phi_{v}(B))\ra
{\bq}/{\bz}.
$$
See \cite{McC}, Theorem 4.8. Thus, since the Pontryagin dual of a
direct sum is a direct product, there exists a nondegenerate pairing
\begin{equation}\label{ggp}
\bigoplus_{v\in U }\Phi_{v}(A)(k(v)) \times\prod_{v\in U }H^{\e
1}\lbe(k(v),\Phi_{v}(B))\ra{\bq}/{\bz}.
\end{equation}
Note that, since the groups $\Phi_{v}(A)(k(v))$ are zero for all but
finitely primes $v$, the above pairing is in fact a pairing of
finite groups. Now, given a nondegenerate pairing of abelian groups
$M\times N\ra {\bq}/{\bz}$ and a subgroup $M^{\prime}\subset M$,
there exists an induced nondegenerate pairing $(M/M^{\prime})\times
(M^{\prime})^{\bot}\ra {\bq}/{\bz}$, where
$(M^{\prime})^{\bot}\subset N$ is the exact annihilator (i.e.,
orthogonal complement) of $M^{\prime}$ in $N$. In the case at hand,
a natural question is the following one: what is the exact
annihilator, under \eqref{ggp}, of the image of $A(F)$ in
$\bigoplus_{v\in U }\Phi_{v}(A)(k(v))$? The main result of this
paper is the determination of this annihilator under the assumption
that $\!\!\Sha^{1}(A)$ is finite. More precisely, the following
holds.

\begin{theorem} Assume that $\!\!\Sha^{1}(A)$ is finite.
Then \eqref{ggp} induces a nondegenerate pairing of finite groups
$$
C_{A,\e F,\e S}\times C_{B,\e F,\e S}^{\e 1}\ra\bq/\bz,
$$
where $C_{B,\e F,\e S}^{\e 1}$ is the group \eqref{c1}.
\end{theorem}

The subgroup $C_{B,\e F,\e S}^{\e 1}$ of $\prod_{\,v\in U }H^{\e
1}\lbe(k(v),\Phi_{v}(B))$ which appears in the statement of the
theorem admits the following simple description. Consider the finite
set of primes
$$
P=\{v\in U\colon\Phi_{v}(A)(k(v))\neq 0\}
$$
and assume, to avoid trivialities, that $P\neq\emptyset$. Then the
nonzero elements of $C_{B,\e F,\e S}^{\e 1}$ are represented by
principal homogeneous spaces for $B$ over $F$ which have a point
defined over $F_{v}$ for every $v\notin P$ and a point defined over
some unramified extension of $F_{v}$ for each $v\in P$, but no point
defined over $F_{v}$ itself for some $v\in P$. Thus the theorem
could be interpreted as saying that the N\'eron class group $C_{A,\e
F,\e S}$ of $A$ controls via duality the existence (or lack of
existence) of rational points on principal homogeneous spaces for
$B$ over $F$ in various completions of $F$.

Regarding the question of how the groups $C_{\lbe A,\e F,\e S}$ and
$\!\!\Sha^{1}(A)$ are related, we show in Proposition 3.4 that there
exists an exact sequence
$$
0\ra C_{A,F,S}\ra D^{\e 1}(U,\sao)\ra
\Sha^{1}\lbe(\lbe A)\ra 0,
$$
where $D^{\e 1}(U,\sao)$ is the group \eqref{Di}. Hence the
Tate-Shafarevich group of $A$ is not only not the correct analogue
of the ideal class group of $F$ but, in a sense, is ``perpendicular" to it.

\begin{remark} The assumption $S\neq\emptyset$ is only used in Section 3, where it guarantees the finiteness of the $S$-class group of an {\it affine} group scheme $H$. However, when $H=A$ is an abelian variety, this hypothesis is not needed and, consequently, our main theorem above remains valid when $S=\emptyset$ in the function field case.
\end{remark}

\section*{Acknowledgements}
I thank B.Conrad, B.Edixhoven and X.Xarles for helpful
correspondence. I also thank the referee for the many suggestions
that led to an improvement of the original presentation.

\section{Preliminaries}

Let $F$ be a global field, i.e. $F$ is a finite extension of
${\mathbb Q}$ (the number field case) or is finitely generated and
of transcendence degree 1 over a finite field of constants (the
function field case). We fix a separable algebraic closure $\fb$ of
$F$. Further, $S$ is a nonempty finite set of primes of $F$
containing the archimedean primes in the number field case, $\s
O_{F,\e S}$ is the corresponding ring of $S$-integers of $F$ and
$U=\spec\s O_{F,\e S}$. For every $v\in U$, ${\s O}_{v}$ will denote
the {\it completion} of the local ring of $U$ at $v$, $F_{v}$ will
denote its field of fractions and $k(v)$ is the corresponding
(finite) residue field. Further, for each prime $v$ of $F$, we fix a
prime $\overline{v}$ of $\fb$ lying above $v$, let
$\fb_{\be\overline{v}}$ be the completion of $\fb$ at $\overline{v}$
(which is a separable algebraic closure of $F_{v}$) and set $G_{\be
F_{v}}=\text{Gal}\big(\e\fb_{\be\overline{v}}/F_{v})$. The inertia
subgroup of $G_{\be F_{v}}$ is the group
$I_{\e\overline{v}}=\text{Gal}\big(\e\fb_{\be\overline{v}}/
F_{v}^{\e\text{nr}}\e\big)$, where $F_{v}^{\e\text{nr}}$ is the
maximal unramified extension of $F_{v}$ lying inside
$\fb_{\be\overline{v}}$. We will write $i_{v}$ for the canonical
closed immersion $\spec k(v)\ra U$.

\smallskip

For any abelian group $M$ and positive integer $n$, we will write
$M_{n}$ for the $n$-torsion subgroup of $M$ and $M/n$ for the
quotient $M/nM$. Further, $M(n)=\cup_{\e r\geq 1}M_{n^{\le r}}$ is
the $n$-primary torsion component of $M$ and
$M_{\text{div}}=\bigcap_{\, n}nM$ is the subgroup of $M$ of
infinitely divisible elements. For simplicity, we will write
$M/\text{div}$ for $M/M_{\text{div}}$. Further, we define
$M^{D}=\text{Hom}(M,\bq/\e\bz)$.

A pairing of abelian groups $A\times B\ra\bq/\e\bz$ is called {\it
nondegenerate on the right} (resp. {\it left}) if the induced
homomorphism $B\ra A^{D}$ (resp. $A\ra B^{D}$) is injective. It is
called {\it nondegenerate} if it is nondegenerate both on the
right and on the left. The pairing is said to be {\it perfect} if
the induced homomorphisms $B\ra A^{D}$ and $A\ra B^{D}$ are
isomorphisms. Note that a pairing of finite abelian groups is
nondegenerate if and only if it is perfect.

For the convenience of the reader, we now recall the well-known
4-lemmas from Homological Algebra.

\begin{lemma} Let
$$
\xymatrix{X \ar[r] \ar@{->>}[d]^(.45){\alpha} & Y \ar[r]
\ar@{-->>}[d]_(.45){\beta} & Z \ar[r] \ar@{->>}[d]^(.45){\gamma}&
\,W\ar@{^{(}->}[d]^{\delta} \\
X^{\prime} \ar[r]& Y^{\prime} \ar[r] & Z^{\prime} \ar[r] &
W^{\prime}}
$$
be an exact commutative diagram of abelian groups and group
homomorphisms. If the maps $\alpha$ and $\gamma$ are epimorphisms
and $\delta$ is a monomorphism, then the map $\beta$ is an
epimorphism. Dually, if
$$
\xymatrix{X \ar[r] \ar@{->>}[d]^(.45){\alpha} & Y \ar[r]
\ar@{^{(}->}[d]_(.45){\beta} & Z \ar[r]
\ar@{^{(}-->}[d]^(.45){\gamma}&
\,W\ar@{^{(}->}[d]^{\delta} \\
X^{\prime} \ar[r]& Y^{\prime} \ar[r] & Z^{\prime} \ar[r] &
W^{\prime}}
$$
is an exact commutative diagram in which $\alpha$ is an epimorphism
and $\beta$ and $\delta$ are monomorphisms, then $\gamma$ is a
monomorphism.\qed
\end{lemma}

\begin{lemma}
Let $n$ be a positive integer and let
$$
A_{1}\overset{f_{1}}\longrightarrow A_{2}\overset{f_{2}}
\longrightarrow A_{3}\overset{f_{3}}\longrightarrow
A_{4}\overset{f_{4}}\longrightarrow  A_{5}
$$
and
$$
B_{5}\overset{g_{4}}\longrightarrow B_{4}\overset{g_{3}}
\longrightarrow B_{3}\overset{g_{2}}\longrightarrow
B_{2}\overset{g_{1}}\longrightarrow  B_{1}
$$
be exact sequences of abelian groups. Assume that there exist
pairings
$$
\varphi_{i}\colon A_i \times B_i \ra \bq/\bz\quad(1\leq i\leq 5)
$$
such that the following conditions hold.
\begin{enumerate}
\item[(i)] The map $A_1/n \ra [(B_1)_n]^{D}$ induced by
$\varphi_1$ is surjective.
\item[(ii)] The maps $A_i/n \ra [(B_i)_n]^{D}$ and $B_i/n\ra
[(A_i)_n]^{D}$ induced by $\varphi_{i}$ are injective for $i=2$ and $i=4$.
\item[(iii)] The map $(A_5)_n \ra [B_5/n]^{D}$ induced by
$\varphi_5$ is injective.
\end{enumerate}
Then the maps $A_3/n \ra [(B_3)_n]^{D}$ and $B_3/n \ra
[(A_3)_n]^{D}$ induced by $\varphi_3$ are injective.
\end{lemma}
\begin{proof} Hypothesis (i) implies that the canonical map $d\colon\img f_1/n
\ra[(\img g_1)_n]^{D}$ is surjective. Now Lemma 2.1 applied to the
exact commutative diagram
$$
\xymatrix{\img f_1/n \ar[r] \ar@{->>}[d]^(.45){d} & A_2/n \ar[r]
\ar@{^{(}->}[d]_(.45){\text{(ii)}} & \img f_2/n \ar[r]
\ar@{^{(}-->}[d]^(.45){b}&
\,0\ar@{^{(}->}[d]\\
[(\img g_1)_n]^D \ar[r] & [(B_2)_n]^D \ar[r] & [(\img g_2)_n]^D
\ar[r] &0}
$$
(whose top and bottom rows are induced by the short exact sequences
$0\ra\img f_1 \ra A_2 \ra\img f_2\ra 0$ and $0\ra\img g_2 \ra B_2
\ra\img g_1 \ra 0$, respectively) shows that the map labeled $b$
above is injective. On the other hand, applying Lemma 2.1 to the
diagrams
$$
\xymatrix{ 0 \ar[r] \ar@{->>}[d]& (\img f_3)_n \ar[r]
\ar@{-->>}[d]^(.45){a} & (A_4)_n \ar[r]
\ar@{->>}[d]_(.45){\text{(ii)}} &
(\img f_4)_n \ar@{^{(}->}[d]_(.45){\text{(iii)}} \\
0 \ar[r] & (\img g_3/n)^{D} \ar[r] & (B_4/n)^{D} \ar[r] & (\img
g_4/n)^{D} }
$$
and
$$
\xymatrix{ (A_4)_n \ar[r] \ar@{->>}[d] & (\img f_4)_n \ar[r]
\ar@{^{(}->}[d] & \img f_3/n \ar[r] \ar@{^{(}-->}[d]^(.45){c}&
A_4/n \ar@{^{(}->}[d]_(.45){\text{(iii)}} \\
[(B_4)_n]^D \ar[r] & (\img g_4/n)^D \ar[r] & [(\img g_3)_n]^D \ar[r]
&
[(B_4)_n]^{D} \\
}
$$
shows that the maps labeled $a$ and $c$ above are surjective and
injective, respectively. Further, applying Lemma 2.1 to the diagram
$$
\xymatrix{ (\img f_3)_n \ar[r] \ar@{->>}[d]^(.45){a} & \img f_2/n
\ar[r] \ar@{^{(}->}[d]^(.45){b} & A_3/n \ar[r] \ar[d]&
\img f_3/n \ar@{^{(}->}[d]^(.45){c} \\
(\img g_3/_n)^{D} \ar[r] & [(\img g_2)_n]^{D} \ar[r] & [(B_3)_n]^{D}
\ar[r] &
[(\img g_3)_n]^{D}, \\
}
$$
we obtain the injectivity of $A_3/n \ra [(B_3)_n]^{D}$. To check the
injectivity of $B_3/n \ra [(A_3)_n]^{D}$, i.e., the surjectivity of
$(A_3)_n\ra [B_3/n]^{D}$, one considers the diagrams
$$
\xymatrix{ (A_2)_n \ar[r] \ar@{->>}[d]_(.45){\text{(ii)}} & (\img f_2)_n
\ar[r] \ar@{-->>}[d]& \img f_1/n \ar[r] \ar@{->>}[d]^(.45){d}&
A_2/n \ar@{^{(}->}[d]_(.45){\text{(ii)}}  \\
(B_2/n)^{D} \ar[r] & (\img g_2/n)^{D} \ar[r] & [(\img g_1)_n]^{D}
\ar[r] & [(B_2)_n]^{D} }
$$
and
$$
\xymatrix{ (\img f_2)_n \ar[r] \ar@{->>}[d] & (A_3)_n \ar[r]
\ar@{-->>}[d]& (\img f_3)_n \ar[r] \ar@{->>}[d]^(.45){a}&
\img f_2/n \ar@{^{(}->}[d]^(.45){b} \\
(\img g_2/n)^{D} \ar[r] & (B_3/n)^{D} \ar[r] & (\img g_3/n)^{D}
\ar[r] & [(\img g_2)_n]^{D}.}
$$
\end{proof}

\section{Class groups and Tate-Shafarevich groups}

In the remainder of the paper we will need to consider flat
cohomology groups $H^{\e i}_{\text{fl}}(U,\s F)=H^{\e
i}(U_{\text{fl}},\s F)$, where $U_{\text{fl}}$ is the category of
$U$-schemes locally of finite type endowed with the flat topology.
If $\s F$ is represented by a smooth, quasi-projective and
commutative $U$-group scheme, then $H^{\e i}_{\text{fl}}(U,\s
F)=H^{i}_{\text{\'et}}(U,\s F)$ (see \cite{Mi1}, Theorem III.3.9,
p.114). On the other hand, if $Y$ is any scheme, $Y_{0}$ will denote
the set of closed points of $Y$. Further, if $V$ is a scheme, $Y$ is
a $V$-scheme and $\spec A$ is an affine $V$-scheme, $Y(A)$ will
denote $\text{Hom}_{_V}\be(\spec A,Y)$.

Let $V$ be a nonempty open subscheme of $U$.  The ring of
$V$-integral adeles of $U$ is by definition
$$
{\mathbb A}_{\e U}(V)=\displaystyle\prod_{v\e\in\e U\e\setminus
V}F_{v}\times\displaystyle\prod_{v\e\in\e V_{0}}{\s O}_{v}.
$$
Since $V=\spec\!\left(\e\bigcap_{\e v\e\in\e V_{0}}{\s
O}_{v}\right)$ and there exists a canonical map $\bigcap_{\e
v\e\in\e V_{0}}{\s O}_{v}\ra{\mathbb A}_{\e U}(V)$, $\spec {\mathbb
A}_{\e U}(V)$ is canonically an affine $V$-scheme. Now let $H$ be a
smooth algebraic group over $F$ and let $\s M$ be a quasi-projective
$U$-model of $H$ of {\it finite type}. We will write $\s M_{_V}\be$
for $\s M\times_{U}V$. The projections ${\mathbb A}_{\e U}(V)\ra
F_{v}$ (for $v\in\e U\setminus V$) and ${\mathbb A}_{\e U}(V)\ra {\s
O}_{v}$ (for $v\e\in\e V_{0}$) induce a bijection
$$
\s M_{_V}({\mathbb A}_{\e U}(V))\ra \displaystyle\prod_{v\e\in\e U\e\setminus
V}H(F_{v})\times\displaystyle\prod_{v\e\in\e V_{0}}\s M_{v}({\s O}_{v}),
$$
where $\s M_{v}=\s M\times_{U}\spec\sov$. See \cite{C2}, Theorem 3.5.
We now define a partial ordering on the family of nonempty open
subschemes of $U$ by setting $V\leq V^{\e\prime}$ if
$V^{\e\prime}\subset V$. Then, for every pair $V,V^{\e\prime}$ of
such schemes such that $V\leq V^{\e\prime}$,
there exists a canonical injection ${\mathbb A}_{\e
U}(V)\hookrightarrow {\mathbb A}_{\e U}(V^{\e\prime})$, namely the
product of the identity map on $\prod_{\,v\in U\e\setminus V}F_{v}$
and the canonical injection
$$
\displaystyle\prod_{\,v\in
V_{0}}\sov\hookrightarrow\displaystyle\prod_{v\e\in\e V\setminus
V^{\e\prime}}F_{v}\e\times\displaystyle\prod_{v\e\in\e
V^{\e\prime}_{0}}\sov.
$$
We will view ${\mathbb A}_{\e U}(V)$ as a subring of ${\mathbb
A}_{\e U}(V^{\e\prime})$ through the above map. The {\it ring of
adeles of $\,U$} is by definition
$$
{\mathbb A}_{\e U}=\varinjlim_{V}\e{\mathbb A}_{\e
U}(V).
$$
Clearly, for any $V$ as above, we may regard $\e{\mathbb A}_{\e
U}(V)$ as a subring of ${\mathbb A}_{\e U}$ and $\s M_{_V}({\mathbb
A}_{\e U}(V))$ as a subgroup of $\s M({\mathbb A}_{\e U})$. Then, by
\cite{C2}, p.5, the natural map
$$
\varinjlim_{V}\,\s M_{_V}({\mathbb A}_{\e U}(V))\ra \s M({\mathbb A}_{\e U})
$$
is a bijection. We conclude that there exists a canonical bijection
\begin{equation}\label{1}
\s M({\mathbb A}_{\e U})=\varinjlim_{V}\left(
\displaystyle\prod_{v\e\in\e U\e\setminus
V}H(F_{v})\times\displaystyle\prod_{v\e\in\e V_{0}}\s M_{v}({\s
O}_{v})\right).
\end{equation}

\begin{remark} For each non-archimedean prime $v$ of $F$, $H(F_{v})$
has a natural locally compact Hausdorff topology containing $\s
M_{v}({\s O}_{v})$ as a compact open subgroup. Thus,
$\varinjlim_{\,V}\be\prod_{\e v\e\in\e U\setminus
V}H(F_{v})\e\times\e\prod_{\e v\e\in\e V_{0}}\s M_{v}({\s O}_{v})$
has a natural locally compact Hausdorff topology \cite{HR}, 6.16(c),
p.57. This topology can then be transferred to $\s M({\mathbb
A}_{U})$ via \eqref{1} so that \eqref{1} is a homeomorphism. See
\cite{C2}, Theorem 3.5.
\end{remark}

It is shown in \cite{C2}, Theorem 4.3, that $H(F)$ injects into $\s
M({\mathbb A}_{\e U})$. Further, as noted above, $\s M({\mathbb
A}_{\e U}(U))\subset\s M({\mathbb A}_{\e U})$. We define the {\it
class set} $\e C(\s M)$ of $\s M$ as the double coset space
$$
C(\s M)=\s M({\mathbb A}_{\e U}(U))\e\backslash\e\s
M\!\left(\lbe{\mathbb A}_{\e U}\lbe\right)\be/\e H(\lbe F).
$$

Now, although the arguments of \cite{Nis1}, Chapter I, \S2 (which
are reproduced in \cite{GA2}, \S3) are valid in principle only when
$\s M$ is affine, they admit a straightforward generalization to
arbitrary $\s M$ as above\footnote{Note, however, that if
topological considerations are relevant, then this generalization
is nontrivial. See \cite{C2}, \S3, and Remark 3.1 above.}. In
particular, the pointed set $C(\s M)$ admits the following
Nisnevich-cohomological interpretation (see \cite{GA2}, Theorem
3.5):
$$
C(\s M)=H_{{\rm{Nis}}}^{1}(U,{\s M}).
$$

Assume now that, in addition to being smooth, $H$ is {\it
commutative, connected and admits a N\'eron model\,}\footnote{In
\cite{BLR}, Chapter 10, this is called a {\it N\'eron lft-model}.}
$\s H$ over $U$. Thus $\s H$ is a smooth and separated $U$-group
scheme (in particular, it is locally of finite type) and represents
the sheaf $j_{\e *}H$ on the small smooth site over $U$. Its
identity component $\sho$ is a smooth $U$-model of $H$ of {\it
finite type} and {\it quasi-projective}. See \cite{SGA3},
$\text{VI}_{\text{B}}$, Proposition 3.9, p.344, and \cite{BLR},
Theorem 6.4.1, p.153. We call the corresponding class group
$C(\sho)$ the {\it N\'eron $S$-class group of $H$} and denote it by
$C_{H,\e F,\e S}$. Thus
$$
C_{H,\e F,\e S}=\sho\!\lbe\left(\lbe{\mathbb A}_{\e
U}\lbe\right)\be/\e H(\lbe F\lbe)\e\sho\!\lbe\left(\lbe{\mathbb A}_{\e U}(U)
\lbe\right).
$$
The group $C_{H,\e F,\e S}$ is known to be finite if $H$ is affine
(see, e.g., \cite{C1}, \S1.3). It is also finite if $H=A$ is an
abelian variety (this is immediate from Theorem 3.2 below since
$\oplus_{\,v\in U }\,\Phi_{v}(A)(k(v))$ is finite).

For each prime $v\in U $, let $\Phi_{v}(H)=i_{v}^{*}\big(\le\s
H/\sho\big)$ be the $k(v)$-sheaf of connected components of $\s H$
at $v$. It is representable by an \'etale $k(v)$-group scheme of
finite type, and there exists a canonical exact sequence of \'etale
sheaves on $U$
\begin{equation}\label{ses}
0\ra\sho\ra\s H\ra\bigoplus_{v\in U } \,\,(i_{v})_{*}\Phi_{v}(H)\ra
0.
\end{equation}
For each $v\in U $, the natural map $\s
H_{\sov}\ra(i_{v})_{*}\Phi_{v}(H)$ induces a map $\vartheta_{\lbe
v}=\vartheta_{\lbe H,\e v}\colon H(F_{v})\ra\Phi_{v}(H)(k(v))$. Let
\begin{equation}\label{2}
\vartheta_{\!_{H,\e S}}\colon H(F)\ra\bigoplus_{v\in U
}\Phi_{v}(H)(k(v))
\end{equation}
be the map induced by the composition
$$
H(F)\ra\prod_{v\in U }H(F_{v})\overset{\prod\vartheta_{\be
v}}\longrightarrow\prod_{v\in U }\Phi_{v}(H)(k(v)),
$$
where the first map is the natural ``diagonal
homomorphism"\footnote{That the image of $\vartheta_{\!_{H,S}}$ is
contained in $\bigoplus_{v\in U }\Phi_{v}(H)(k(v))$ follows from
\cite{T}, Lemma 10.1.1, p.157.}.

\begin{teorema} Let $H$ be a smooth, connected and commutative algebraic
group over $F$. Assume that $H$ admits a N\'eron model $\s H$ over
$U$. Then \eqref{ses} induces an exact sequence
$$
0\ra\sho(U)\ra H(F)\overset{\vartheta_{\!_{H\be,S}}}\longrightarrow
\displaystyle\bigoplus_{v\in U }\,\Phi_{v}(H)(k(v))\ra  C_{\e H,\e
F,\e S}\ra 0,
$$
where $\vartheta_{\!_{H,S}}$ is the map \eqref{2}.
\end{teorema}
\begin{proof} Taking \'etale cohomology of \eqref{ses} and using the identification
$$
C_{H,\e F,\e S}=H_{{\rm{Nis}}}^{1}(U,\sho),
$$
we are immediately reduced to checking that
$H_{{\rm{Nis}}}^{1}(U,\sho)$ is canonically isomorphic to the kernel
of the natural map $H^{1}_{{\rm{\acute{e}t}}}\big(\e U,\sho\big)\ra
H^{1}_{{\rm{\acute{e}t}}}\big(\e U,\s H)$. The Cartan-Leray spectral
sequence shows that the maps $H^{1}_{{\rm{\acute{e}t}}}\big(\e
U,\sho\big)\ra H^{1}_{{\rm{\acute{e}t}}}\big(\e U,\s H)$ and $
H^{1}_{{\rm{\acute{e}t}}}\be\big(\e U,\sho\big)\ra
H^{1}_{{\rm{\acute{e}t}}}\be\big(\e F,H)$ have the same kernel (see
\cite{GA1}, proof of Lemma 3.1). On the other hand, Ye.Nisnevich has
shown \cite{Nis2}, Example 1.44, p.286, that the canonical sequence
$$
0\ra H_{{\rm{Nis}}}^{1}(U,\sho)\ra H^{1}_{{\rm{\acute{e}t}}}\big(\e
U,\sho\big)\ra H^{1}(F,H\e)
$$
is exact provided the local maps $H^{1}_{{\rm{\acute{e}t}}}(\sov^{\e
h},\sho)\ra H^{1}_{{\rm{\acute{e}t}}}(F_{v}^{\e h},\sho)$ are
injective for each $v\in U_{0}$, where $\sov^{\e h}$ denotes the
henselization of the local ring of $U$ at $v$ and $F_{v}^{\e h}$ is
its field of fractions. Since $H^{1}_{{\rm{\acute{e}t}}}(\sov^{\e
h},\sho)$ is in fact zero \cite{Mi1}, Remark III.3.11, p.116, and
Lang's Theorem \cite{L}, the proof is now complete.
\end{proof}

\begin{remark} The literature records the following assertions:
``If $A$ is an abelian variety over a number field $F$, then $A(F)$
is the natural analog of the group of units of $F$ and the
Tate-Shafarevich group of $A$ is the natural analog of the ideal
class group of $F$". We believe that these assertions are incorrect.
Indeed, the exact sequence of the theorem for $H=\bg_{m,F}$ is the
familiar exact sequence
\begin{equation}\label{sos}
1\ra \s O_{F,S}^{*}\ra F^{*}\ra\displaystyle\bigoplus_{v\in U
}\,\mathbb{Z}\ra C_{F,S}\ra 0,
\end{equation}
where $\sho(U)=\bg_{m,U}(U)=\s O_{F,S}^{*}$ is the group of
$S$-units of $F$, $\bigoplus_{v\in U }\mathbb{Z}$ is (isomorphic to)
the group $\s I_{F,S}$ of fractional $S$-ideals of $F$ and $C_{F,S}$
is the ideal $S$-class group of $F$. See \cite{GA1}, Remark 2.1.
Thus \eqref{sos} and the exact sequence of the theorem when $H=A$ is
an abelian variety with N\'eron model $\s A$ show that $\sao(U)$,
$\oplus_{v\in U }\Phi_{v}(A)(k(v))$ and $C_{\e A,\e F,\e S}$ are
natural analogs of $\s O_{F,S}^{*}$, $\s I_{F,S}$ and $C_{F,S}$,
respectively.
\end{remark}

\smallskip

Let $H$ be as in the statement of the theorem. For each $v \notin S$, set
$$
H^{\e 1}_{\text{nr}}(F_{v},H)=\krn\!\!\left[H^{\e 1}(F_{v},H)\ra H^{\e 1}
(F^{\text{nr}}_{v},H)\right],
$$
where the map involved is the restriction map in Galois cohomology. The inflation map in Galois cohomology induces an isomorphism
$$
H^{\e 1}\lbe(G_{\lbe
F_{v}}\be/I_{\overline{v}},H(F^{\text{nr}}_{v}))\simeq H^{\e
1}_{\text{nr}}(F_{v},H).
$$
On the other hand, by a straightforward generalization of
\cite{Mi2}, proof of Proposition I.3.8, p.57, the reduction map
$H(F^{\e\text{nr}}_{v})\ra\Phi_{v}(H)\big(\e\overline{k(v)}\e \big)$
induces an isomorphism $H^{\e 1}(G_{\be
F_{v}}\be/I_{\overline{v}},H(F^{\text{nr}}_{v}))\simeq H^{\e
1}\be(k(v),\Phi_{v}(H))$. Thus there exists a canonical isomorphism
\begin{equation}\label{Finr}
H^{\e 1}_{\text{nr}}(F_{v},H)\simeq H^{\e 1}\be(k(v),\Phi_{v}(H)).
\end{equation}
We will henceforth identify $H^{\e 1}_{\text{nr}}(F_{v},H)$ and
$H^{\e 1}\be(k(v),\Phi_{v}(H))$ via the above map.

Now, for {\it any} prime $v$ of $F$, there exists a canonical
``localization" map $H^{\e 1}(F,H)\ra H^{\e 1}(F_{v},H)$. We let
\begin{equation}\label{loc1}
\lambda_{S} \colon H^{\e 1}(F,H)\ra\displaystyle\prod_{v\in U }H^{\e
1}(F_{v},H)
\end{equation}
be the induced map and set $\!\!\Sha^{1}_{S}(H)=\krn\lambda_{S}$.
Further, we recall the Tate-Shafarevich group of $H$:
$$
\Sha^{1}(H)=\krn\!\left[\e H^{\e
1}(F,H)\ra\displaystyle\prod_{\text{all $v$}}H^{\e
1}(F_{v},H)\right].
$$
Now define
$$
\Sha_{S^{\lbe\text{c}}}^{1}(H)=\krn\!\!\left[ H^1(F,H) \ra
\displaystyle\prod_{v \in S}H^{\e 1}(F_{v},H)\right]
$$
and let
$$
\lambda_{S}^{\e\prime}\colon\Sha_{S^{\lbe\text{c}}}^{1}(H)\ra
\prod_{v\in U }H^1(F_v,H)
$$
be the restriction of \eqref{loc1} to $\!\!\Sha_{S^{\lbe\text{c}}}^{1}(H)$. Set
\begin{equation}\label{c1}
C_{H,F,S}^{\e 1}=\left(\e\prod_{v \notin S}
H_{\text{nr}}^{1}(F_{v},H)\right)\cap \img \lambda_{S}^{\prime}
\subset \prod_{v \notin S}
H^{1}(F_{v},H).
\end{equation}
The elements of $C_{H,F,S}^{\e 1}$ can be described as follows.
Consider the set of primes
$$
P=\{v\in U \colon H^{\e 1}\lbe(k(v),\Phi_{v}(H))\neq 0\}
$$
and assume, to avoid trivialities, that $P\neq\emptyset$. Then the
nonzero elements of $C_{H,F,S}^{\e 1}$ are represented by principal
homogeneous spaces for $H$ over $F$ which have a point defined over
$F_{v}$ for every $v\notin P$ and a point defined over some
unramified extension of $F_{v}$ for each $v\in P$, but no point
defined over $F_{v}$ itself for some $v\in P$.

\smallskip

In order to explain how the groups $C_{H,F,S}$, $C_{H,F,S}^{\e 1}$
and $\!\!\Sha^{1}\be(\lbe H)$ are related, we need another
definition. For any prime $v$ of $F$ and any sheaf $\s F$ on
$U_{\text{fl}}$, let $H^{\e i}_{\text{fl}}(F_{v},\s F)=H^{\e
i}_{\text{fl}}(\spec F_{v},\s F_{v})$, where $\s F_{v}$ is the sheaf
on $(\spec F_{v})_{\text{fl}}$ obtained by pulling back $\s F$
relative to the composite morphism
$$
\spec F_{v}\ra\spec F\ra U.
$$
If $v$ is archimedean, $H^{\e i}_{\text{fl}}(\spec F_{v},\s F_{v})$
will denote the $i$-th {\it reduced} (Tate) cohomology group of $\s
F_{v}$. Clearly, the preceding morphism induces a map $H^{\e
i}_{\text{fl}}(U,\s F)\ra H^{\e i}_{\text{fl}}(F_{v},\s F)$. Set
\begin{equation}\label{Di}
D^{\e i}(U,\s F)=\krn\!\!\left[H^{\e i}_{\text{fl}}(U,\s
F)\ra\prod_{v\in S} H^{\e i}_{\text{fl}}(F_{v},\s F)\right].
\end{equation}
Recall that $H^{\e i}_{\text{fl}}(U,\s F)=H^{\e
i}_{\text{\'et}}(U,\s F)$ if $\s F$ is represented by a smooth,
quasi-projective and commutative $U$-group scheme.

\begin{proposition} Let $H$ be a smooth, connected and commutative algebraic
group over $F$. Assume that $H$ admits a N\'eron model $\s H$ over
$U$. Then there exist canonical exact sequences
$$
0\ra C_{H,F,S}\ra D^{\e 1}(U,\sho)\ra
\Sha^{1}\lbe(\lbe H)\ra 0
$$
and
$$
0\ra \Sha^{1}\lbe(\lbe H)\ra D^{\e 1}(U,\s H)\ra
C_{H,F,S}^{\e 1} \ra 0.
$$
\end{proposition}
\begin{proof} The Cartan-Leray spectral sequence \cite{Mi1}, Theorem 1.18(a), p.89,
and \cite{Mi2}, Remark I.3.10, p.58, yield an exact sequence
$$
0\ra H^{\e 1}_{\text{\'et}}(U,\s H)\ra H^{\e 1}(F,H)\ra \prod_{v\in
U }H^{\e 1}(F^{\text{nr}}_{v},H).
$$
See \cite{Mi2}, proof of Lemma II.5.5, p.247. Now, applying the kernel-cokernel exact sequence (\cite{Mi2}, Proposition I.0.24, p.19) to the pair of maps
$$
H^{\e 1}(F,H)\ra \prod_{v\in U }H^{\e 1}(F_{v},H)\ra\prod_{v\in U
}H^{\e 1}(F^{\text{nr}}_{v},H),
$$
we obtain an exact sequence
$$
0\ra\Sha^{1}_{S}(H)\ra H^{\e 1}(U,\s H)\ra \prod_{v\in U }H^{\e
1}_{\text{nr}}(F_{v},H)\ra\cok\lambda_{S},
$$
where $\lambda_{S}$ is the localization map \eqref{loc1}.
Thus, there exists an exact sequence
\begin{equation}\label{esq1}
0\ra\Sha^{1}_{S}(H)\ra H^{\e 1}(U,\s H)\ra \left(\,\prod_{v\in U
}H^{\e 1}_{\text{nr}}(F_{v},H)\right)\cap\img\lambda_{S}\ra 0.
\end{equation}
On the other hand, via the identification \eqref{Finr}, the exact sequence \eqref{ses} induces an exact sequence
$$
0\ra C_{H,F,S}\ra H^{\e 1}(U,\sho)\ra H^{\e
1}(U,\s H)\ra\prod_{v \notin S}H^{\e
1}_{\text{nr}}(F_{v},H).
$$
Thus \eqref{esq1} yields an exact sequence
\begin{equation}\label{esq2}
0\ra C_{H,F,S}\ra H^{\e 1}(U,\sho)\ra
\Sha^{1}_{S}(\lbe H)\ra 0.
\end{equation}
The proposition now follows from \eqref{c1} and the exact commutative diagrams
$$
\xymatrix{\Sha^{1}_{S}(H)\ar[d]\ar@{^{(}->}[r] &H^{\e 1}(U,\s H)\ar[d]\ar@{->>}[r] &\left(\,\displaystyle\prod_{v \notin S}H^{\e 1}_{\text{nr}}(F_{v},H)\right)\cap\img\lambda_{S}\\
\displaystyle\prod_{v \in S}H^{\e 1}(F_{v},H)\ar@{=}[r] &
\displaystyle\prod_{v \in S}H^{\e 1}(F_{v},H) &&}
$$
and
$$
\xymatrix{0\ar[r] &C_{H,F,S}\ar[r] &H^{\e 1}(U,\sho)\ar[d]\ar[r] &\Sha^{1}_{S}(\lbe H)\ar[d]\ar[r] &0\\
&&\displaystyle\prod_{v \in S}H^{\e 1}(F_{v},H)\ar@{=}[r] &
\displaystyle\prod_{v \in S}H^{\e 1}(F_{v},H),}
$$
whose top rows are \eqref{esq1} and \eqref{esq2}, respectively.
\end{proof}

\section{Proof of the main theorem}

Let $A$ be an abelian variety defined over $F$ and let $B$ be the
abelian variety dual to $A$. The N\'eron model of $B$ over $U$ will
be denoted by $\s B$.

For each $v\in U $, the \'etale $k(v)$-sheaves $\Phi_{v}(A)$ and
$\Phi_{v}(B)$ will be identified with the $G_{k(v)}$-modules
$\Phi_{v}= \Phi_{v}(A)\big(\e\overline{k(v)}\e\big)$ and
$\Phi_{v}^{\e\prime}=\Phi_{v}(B)\big(\e\overline{k(v)}\e\big)$,
respectively. Let $\gv$ and $\gv^{\e\prime}$ be subsheaves of
$\Phi_{v}(A)$ and $\Phi_{v}(B)$. These will be regarded as
$G_{k(v)}$-submodules of $\Phi_{v}$ and $\Phi_{v}^{\e\prime}$,
respectively. We will assume throughout that
\begin{center}
$\gv^{\e\prime}$ is the exact annihilator of $\gv$
\end{center}
under Grothendieck's pairing
$$
\Phi_{v}\times\Phi_{v}^{\e\prime}\ra{\bq}/{\bz}.
$$
Since the above pairing is known to be nondegenerate in the context
of this paper (see, e.g., \cite{McC}, Theorem 4.8),
$\gv^{\e\prime}=\Phi_{v}^{\e\prime}$ and $\gv=0$ is a valid choice.
We will make this choice in the proof of Corollary 4.8 below, and
Corollary 4.8 is used in the proof of Theorem 4.9 (the main
theorem). It is in this way that the nondegeneracy of Grothendieck's
pairing intervenes in the derivation of our main result.

We will need the cohomology groups with compact support $H^{\e
i}_{\text{c}}(U_{\text{fl}},\s F)$ defined in \cite{Mi2},
pp.270-271. For any sheaf $\s F$ on $U_{\text{fl}}$, there exists an
exact sequence
\begin{equation}\label{ccs}
\dots\ra H^{\e i}_{\text{c}}(U_{\text{fl}},\s F\e)\ra H^{\e
i}_{\text{fl}}(U,\s F\e)\ra \displaystyle\prod_{v\in S}H^{\e
i}_{\text{fl}}(F_{v},\s F)\ra H^{\e i+1}_{\text{c}}(U_{\text{fl}},\s
F\e)\ra\dots.
\end{equation}
See \cite{Mi2}, Remark III.0.6(b), p.274.

Let $\g=\bigoplus_{v\in U }\,(i_{v}\lbe)_{*}\gv$ and
$\g^{\e\prime}=\bigoplus_{v\in U }\, (i_{v}\lbe)_{*}\gv^{\e\prime}$,
and define $\s A^{\g}$ and $\s B^{\e\g^{\e\prime}}$ by the exactness
of the sequences
\begin{equation}\label{Gamma}
0\ra\s A^{\circ}\ra\s A^{\g}\ra\g\ra 0
\end{equation}
and
$$
0\ra\s B^{\e\circ}\ra\s B^{\g^{\e\prime}}\ra\g^{\e\prime}\ra 0.
$$
By \cite{SGA7}, Theorem VIII.7.1(b), the canonical Poincar\'e
biextension of $(A,B)$ by $\bg_{\e m,F}$ extends uniquely to a
biextension of $\big(\lbe\s A^{\g}\lbe,\s B^{\e\g^{\e\prime}}\big)$
by $\bg_{m,U}$. This biextension induces a map
$$
\s A^{\g}\otimes^{\mathbf{L}}\s B^{\e\g^{\e\prime}}\ra\bg_{m,U}[1]
$$
in the derived category of the category of smooth sheaves on U, and
this map induces pairings
\begin{equation}\label{pair}
\langle -,-\rangle\colon H^{\e
1}_{\text{c}}\be\big(U_{\rm{fl}},\sa^{\g}\big)\times
H^{1}_{\rm{fl}}\big(U,\s B^{\e\g^{\e\prime}}\big) \ra\bq/\e\bz
\end{equation}
and
\begin{equation}\label{pair'}
H^{\e 1}_{\rm{fl}}\be\big(U,\sa^{\g}\big)\times
H^{1}_{\text{c}}\big(U_{\rm{fl}},\s B^{\e\g^{\e\prime}}\big)
\ra\bq/\e\bz.
\end{equation}
See \cite{Mi2}, comments preceding Theorem III.0.16, and \cite{GH}.

Now recall the group \eqref{Di}:
$$\begin{array}{rcl}
D^{\e i}(U,\s F\e)&=&\!\krn\!\!\left[H^{\e i}_{\rm{fl}}(U,\s
F\e)\ra\displaystyle\prod_{v\in S}
H^{\e i}_{\rm{fl}}(F_{v},\s F)\right]\\\\
&=&{\rm{Im}}\be\left[\,H^{\e i}_{\e\rm{c}}(U_{\rm{fl}},\s F\e)\to
H^{\e i}_{\rm{fl}}(U,\s F\e)\,\right].
\end{array}
$$
There exists a canonical exact commutative diagram
\[
\xymatrix{0\ar[r]& D^{\e 1}\big(U,\s B^{\e\g^{\e\prime}}\big)\ar[r] & H^{1}_{\rm{fl}}
\big(U,\s B^{\e\g^{\e\prime}}\big)\ar[d]\ar[r] &
\displaystyle\prod_{v\in S} H^{\e 1}_{\rm{fl}}(F_{v},B)\ar[d]\\
&& H^{\e 1}_{\text{c}}\big(U_{\rm{fl}},\s A^{\g}\lbe\big)^{\be
D}\ar[r] & \displaystyle\bigoplus_{v\in S} {H}^{\e
0}_{\rm{fl}}(F_{v},A)^{D},}
\]
where the middle vertical map is induced by \eqref{pair} and the
right-hand vertical map is induced by the local pairings \cite{Mi2},
I.3.4, I.3.7 and III.7.8. The horizontal maps come from \eqref{ccs}
(for appropriate choices of $i$ and $\s F$) and the definition of
$D^{\e 1}\big(U,\s B^{\e\g^{\e\prime}}\big)$. It follows that there
exists a well-defined pairing
\begin{equation}\label{pair2}
\{-,-\}\colon D^{\e 1}\big(U,\sa^{\g}\big)\times D^{\e 1}\big(U,\s
B^{\e\g^{\e\prime}}\lbe\big)\ra{\bq}/{\bz}
\end{equation}
given by $\{a,a^{\e\prime}\}=\langle\e
\tilde{a},a^{\e\prime}\e\rangle$, where $a^{\e\prime}\in D^{\e
1}\big(U,\s B^{\e\g^{\e\prime}}\lbe\big)\subset H^{\e
1}_{\rm{fl}}\big(U,\s B^{\e\g^{\e\prime}}\lbe\big)$ and $\tilde{a}$
is a preimage of $a$ in $H^{\e
1}_{\e\rm{c}}\big(U_{\rm{fl}},\sa^{\g}\big)$. We will show (Theorem
4.7) that, if $\!\!\Sha^{1}(A)$ and $\!\!\Sha^{1}(B)$ are
finite\footnote{Recall that $\!\!\Sha^{1}(A)$ is finite if, and only
if, $\!\!\Sha^{1}(B)$ is finite \cite{Mi2}, Remark I.6.14(c),
p.102.}, then \eqref{pair2} is a perfect pairing of finite groups.

\begin{lemma} Let $V$ be a nonempty open subscheme of $U$ such that
$\s A\!\be\mid_{V}$ is an abelian scheme. Let $n$ be any positive
integer. Then the maps
$$
H_{{\rm{c}}}^{1}(V_{\rm{fl}},\s B)/n\ra(H^{1}_{\rm{fl}}(V,\s
A)_n)^{D}
$$
and
$$
H^{1}_{\rm{fl}}(V,\s A)/n\ra(H_{\text{c}}^{1}(V_{\rm{fl}},\s
B)_n)^{D},
$$
induced by \eqref{pair'} over $V$, are injective.
\end{lemma}

\begin{proof}
The lemma is immediate from the commutativity of the diagrams
$$
\xymatrix{
H_{\text{c}}^{1}(V_{\rm{fl}},\s B)/n \ar@{^{(}->}[r] \ar[d]&
H_{\text{c}}^{2}(V_{\rm{fl}},\s B_n) \ar[d]^(.45){\simeq} \\
(H^{1}_{\rm{fl}}(V,\s A)_n)^{D} \ar[r] & H^{1}_{\rm{fl}}(V,\s
A_n)^{D} }
$$
and
$$
\xymatrix{
H^{1}_{\rm{fl}}(V,\s A)/n \ar@{^{(}->}[r] \ar[d]& H^{2}_{\rm{fl}}(V,\s A_n)
\ar[d]^(.45){\simeq} \\
(H_{\text{c}}^{1}(V_{\rm{fl}},\s B)_n)^{D} \ar[r] &
H_{\text{c}}^{1}(V_{\rm{fl}},\s B_n)^{D}, }
$$
where the right-hand vertical maps are isomorphisms by \cite{Mi2},
corollary III.3.2, p.313, and theorem III.8.2, p.361.
\end{proof}

\begin{lemma} Let $V$ and $n$ be as in the previous lemma. Then there exists a canonical exact sequence
$$
0\ra(\!\!\Sha^1(A)_{\rm{div}}\cap\!\!\Sha^1(A)_{n})^{D} \ra
H_{{\rm{c}}}^{2}(V_{\rm{fl}},\s B)_{n}\ra(H^{\e 0}_{\rm{fl}}(V,\s
A)/n)^{D}\ra 0.
$$
In particular, if $\Sha^{1}(A)(n)$ is finite, then the canonical map
$H_{\rm{c}}^{2}(V_{\rm{fl}},\s B)_{n}\ra(H^{\e 0}_{\rm{fl}}(V,\s
A)/n)^{D}$ is an isomorphism.
\end{lemma}
\begin{proof} There exist exact commutative diagrams
$$
\xymatrix{ 0 \ar[r] & H_{\text{c}}^{1}(V_{\rm{fl}},\s B)/n \ar[r]
\ar[d]& H_{\text{c}}^{2}(V_{\rm{fl}},\s B_n) \ar[r]
\ar[d]^(.48){\simeq}& H_{\text{c}}^{2}(V_{\rm{fl}},\s B)_n \ar[r]
\ar[d]&
0 \\
0 \ar[r] & (H^{1}_{\rm{fl}}(V,\s A)_n)^{D} \ar[r] &
H^{1}_{\rm{fl}}(V,\s A_n)^{D} \ar[r] & (H^{\e 0}_{\rm{fl}}(V,\s
A)/n)^{D} \ar[r] & 0 }
$$
and
$$
\xymatrix{ \displaystyle\prod_{v \notin V}H^{0}_{\rm{fl}}(F_{\be
v},B)/n \ar[r] \ar[d]^{\simeq}& H_{\text{c}}^1 (V_{\rm{fl}}, \s B)/n
\ar[r] \ar[d]& D^{1}(V, \s B)/n \ar[r] \ar[d]&
0 \\
\displaystyle\bigoplus_{v \notin V}\,(H^{1}_{\rm{fl}}(F_v,A)_n)^D
\ar[r] & (H^{1}_{\rm{fl}}(V,\s A)_n)^D \ar[r] & (D^{1}(V,\s A)_n)^D
\ar[r] & 0 }
$$
(cf. \cite{Mi2}, p.245). Further, there exist canonical isomorphisms
$D^1(V,\s A)\simeq\!\!\Sha^1(A)$ and $D^1(V,\s
B)\simeq\!\!\Sha^1(B)$ \cite{Mi2}, Lemma II.5.5 p.247. The lemma now
follows from the above diagrams and the existence of a perfect
pairing $\!\!\Sha^1(A)/\rm{div}\times \Sha^1(B)/\rm{div}\ra\bq/\bz$
\cite{HS}, Theorem 4.8, and \cite{GA3}, Theorem 6.6.
\end{proof}

\begin{proposition} Let $n$ be any integer such that $\!\!\Sha^1(A)(n)$ is finite.
Then the maps
$$
H_{{\rm{c}}}^{1}\be\big(U_{\rm{fl}},\s
B^{\e\g^{\e\prime}}\big)/n\ra(H^{1}_{\rm{fl}}(U,\s A^{\g})_n)^{D}
$$
and
$$
H^{1}_{\rm{fl}}(U,\s A^{\g})/n \ra
(H_{{\rm{c}}}^{1}\be\big(U_{\rm{fl}},\s
B^{\e\g^{\e\prime}}\big)_n)^D,
$$
induced by \eqref{pair'}, are injective.
\end{proposition}

\begin{proof} We wish to apply Lemma 2.2.
Let $V$ be a nonempty open subscheme of $U$ such that $\s
A^{\g}\!\be\mid_{V}=\s A\!\be\mid_{V}$ and $\s
B^{\g^\prime}\!\!\be\mid_{V}=\s B\!\be\mid_{V}$ are abelian schemes.
There exist canonical exact sequences of abelian groups
$$
\displaystyle\bigoplus_{v \in U \setminus V}
H^{0}_{\rm{fl}}\be\big(\s O_v,\s B^{\g^\prime}\big) \ra
H_{{\rm{c}}}^{1}(V_{\rm{fl}},\s B) \ra
H_{{\rm{c}}}^{1}\be\big(U_{\rm{fl}},\s B^{\g^\prime}\big) \ra
\displaystyle\bigoplus_{v \in U \setminus V}
H^{1}_{\rm{fl}}\be\big(\s O_v,\s B^{\g^\prime}\big) \ra
H_{{\rm{c}}}^{2}(V_{\rm{fl}},\s B)
$$
and
$$
H^{0}_{\rm{fl}}(V,\s A)\ra \displaystyle\bigoplus_{v \in U \setminus
V} H_{v}^{1}\be\big(\s O_{v,\e{\rm{fl}}},\s A^{\g}\big) \ra
H^{1}\be\big(U_{\rm{fl}},\s A^{\g}\big) \ra H^{1}_{\rm{fl}}(V,\s A)
\ra \displaystyle\bigoplus_{v \in U \setminus V} H_{v}^{2}\be\big(\s
O_{v,\e{\rm{fl}}},\s A^{\g}\big).
$$
See \cite{Mi2}, Proposition III.0.4(c), p.272, and Remark
III.0.6(b), p.275. Further, for each $v \in U$, there exist pairings
$\varphi_{1,v}\colon H^{0}_{\rm{fl}}\big(\s O_v,\s
B^{\g^\prime}\big) \times H_{v}^{2}\big(\s O_{v,\e{\rm{fl}}}, \s
A^{\g}\big) \ra\bq/\bz $ and $\varphi_{4,v}\colon
H^{1}_{\rm{fl}}\big(\s O_v,\s B^{\g^\prime}\big)\times
H_{v}^{1}\big(\s O_{v,\e{\rm{fl}}},\s A^{\g}\big)\ra\bq/\bz$. The
first of these pairings induces a surjection $H^{0}_{\rm{fl}}\big(\s
O_v,\s B^{\g^\prime}\big)/n\ra \big[H_{v}^{2}\big(\s
O_{v,\e{\rm{fl}}}, \s A^{\g}\big)_n\big]^{D}$, and the second one is
a perfect pairing of finite groups. See \cite{Mi2}, Theorem III.2.7,
p.307, and Theorem III.7.13, p.358. On the other hand, by the
previous lemma, the map $H_{{\rm{c}}}^{2}(V_{\rm{fl}},\s B)_n \ra
(H^{0}_{\rm{fl}}(V,\s A)/n)^D$ induced by the pairing
$H_{{\rm{c}}}^{2}(V_{\rm{fl}},\s B)\times H^{\e 0}_{\rm{fl}}(V,\s
A)\ra\bq/\bz$ is an isomorphism. Finally, Lemma 4.1 shows that all
the conditions needed to apply Lemma 2.2 hold, and the proposition
follows.
\end{proof}

Now let $\ell$ be any prime number and let $r$ be a positive
integer. Set
$$
{\rm{Sel}}\big(U,\s B^{\e\g^{\e\prime}}\lbe\big)_{\ell^{
r}}=\krn\!\be\left[H^{\e 1}_{\rm{fl}}\be\big(U,\s
B^{\e\g^{\e\prime}}_{\ell^{r}}\big)\ra\displaystyle\prod_{v\in
S}H^{\e 1}_{\rm{fl}}(F_{v},B)\right],
$$
where the map involved is the composite
$$
H^{\e 1}_{\rm{fl}}\be\big(U,\s
B^{\e\g^{\e\prime}}_{\ell^{r}}\big)\twoheadrightarrow H^{\e
1}_{\rm{fl}}\be\big(U,\s B^{\e\g^{\e\prime}}\be\big)_{
\ell^{r}}\ra\displaystyle\prod_{v\in S}H^{\e 1}_{\rm{fl}}(F_{\lbe
v},B)
$$
The kernel-cokernel exact sequence \cite{Mi2}, Proposition I.0.24,
applied to the above pair of maps yields an exact sequence
\begin{equation}\label{selmer}
0\ra H^{\e 0}_{\rm{fl}}\be\big(U,\s B^{\e\g^{\e\prime}}\big)/\e
\ell^{r} \ra {\rm{Sel}}\big(U,\s
B^{\e\g^{\e\prime}}\lbe\big)_{\ell^{r}}\ra D^{\e 1}\be\big(U,\s
B^{\e\g^{\e\prime}}\lbe\big)_{\ell^{r}}\ra 0.
\end{equation}

Now, by \cite{SGA7}, VIII.2.2.5, the biextension $\big(\lbe\s
A^{\g}\lbe,\s B^{\e\g^{\e\prime}};\bg_{m,\e U}\big)$ induces a map
$$
\s A^{\g}_{\e \ell^{r}}\times\s
B^{\e\g^{\e\prime}}_{\ell^{r}}\ra\bg_{m,\e U}.
$$
The above map induces a (possibly degenerate) pairing
$$
H^{\e 2}_{\text{c}}\be\big(U_{\rm{fl}},\s A^{\g}_{\e
\ell^{r}}\big)\times H^{\e 1}_{\rm{fl}}\be\big(U,\s
B^{\e\g^{\e\prime}}_{\ell^{r}}\big)\ra{\bq}/{\bz},
$$
and we let
\begin{equation}\label{eta}
\eta\colon H^{\e 2}_{\text{c}}\be\big(U_{\rm{fl}},\s A^{\g}_{\e
\ell^{r}}\big)\ra H^{\e 1}_{\rm{fl}}\be\big(U,\s
B^{\e\g^{\e\prime}}_{\ell^{r}}\big)^{D}
\end{equation}
be the map induced by this pairing.

On the other hand, the exact sequences of flat sheaves $0\ra \s
A^{\g}_{\e \ell^{r}}\ra\s A^{\g}\overset{\ell^{r}}\longrightarrow \s
A^{\g}\ra 0$ and $0\ra \s B^{\e\g^{\e\prime}}_{\e \ell^{r}}\ra\s
B^{\e\g^{\e\prime}}\overset{\ell^{r}}\longrightarrow \s
B^{\e\g^{\e\prime}}\ra 0$ induce maps
\begin{equation}\label{partialc}
\partial_{\e\text{c}}\colon H^{\e 1}_{\text{c}}\be\big(U_{\rm{fl}},\s A^{\g}\big)\ra H^{\e 2}_{\text{c}}(U_{\rm{fl}},\s A^{\g}_{ \ell^{r}}).
\end{equation}
and
\begin{equation}\label{vartheta}
\vartheta \colon H^{\e 1}_{\rm{fl}}\be\big(U,\s
B^{\e\g^{\e\prime}}_{\ell^{r}}\big)\ra H^{\e
1}_{\rm{fl}}\be\big(U,\s B^{\e\g^{\e\prime}}\big)_{\ell^{r}}
\end{equation}
such that the following holds. If
\begin{equation}\label{ev}
[-,-]\colon H^{\e 1}_{\rm{fl}}\be\big(U,\s
B^{\e\g^{\e\prime}}_{\ell^{r}}\big)^{D}\times H^{\e
1}_{\rm{fl}}\be\big(U,\s
B^{\e\g^{\e\prime}}_{\ell^{r}}\big)\ra{\bq}/{\bz}
\end{equation}
is the evaluation pairing and $\langle-,-\rangle$ is the pairing \eqref{pair}, then
\begin{equation}\label{compatible}
[\e\eta\partial_{\text{c}}(\zeta),\xi\e]=\langle\zeta,
\vartheta(\xi)\rangle
\end{equation}
for every $\zeta\in H^{\e 1}_{\text{c}}\be\big(U_{\rm{fl}},\s
A^{\g}\big)$ and $\xi\in H^{\e 1}_{\rm{fl}}\be\big(U,\s
B^{\e\g^{\e\prime}}_{\ell^{r}}\big)$, where
$\eta,\partial_{\e\text{c}}$ and $\vartheta$ are the maps
\eqref{eta}, \eqref{partialc} and \eqref{vartheta}. Now consider
\begin{equation}\label{delta'}
\delta^{\e\prime}\colon\displaystyle{\prod_{v\in S}} H^{\e
1}_{\rm{fl}}(F_{v},A_{\e \ell^{r}}\!)\ra H^{\e
2}_{\text{c}}\be\big(U_{\rm{fl}},\s A^{\g}_{\e \ell^{r}}\big)
\end{equation}
(see \eqref{ccs}) and set
\begin{equation}\label{delta}
\delta=\eta\circ\delta^{\e\prime}\colon \displaystyle{\prod_{v\in
S}} H^{\e 1}_{\rm{fl}}(F_{v},A_{\e \ell^{r}}\!) \ra H^{\e
1}_{\rm{fl}}\be\big(U,\s B^{\e\g^{\e\prime}}_{\ell^{r}}\big)^{\lbe
D}.
\end{equation}
By \cite{Mi2}, Corollary I.2.3, p.34, and Theorem III.6.10, p.344,
there exists a perfect ``cup-product" pairing\footnote{Recall that,
in general, direct products and direct sums are in natural duality.
Note, however, that $\prod_{\e v\in S} H^{\e
1}_{\rm{fl}}(F_{v},A_{\e \ell^{r}}\!)$ and $\bigoplus_{\e v\in S}
H^{\e 1}_{\rm{fl}}(F_{v},A_{\e \ell^{r}}\!)$ are canonically
isomorphic since $S$ is finite.}
$$
(-,-)\,\colon\displaystyle{\prod_{v\in S}}\,H^{\e
1}_{\rm{fl}}(F_{v},A_{
\ell^{r}}\!)\times\displaystyle{\bigoplus_{v\in S}}\,H^{\e
1}_{\rm{fl}}(F_{v},B_{ \ell^{r}}\!)\ra{\bq}/{\bz}
$$
defined by
$$
\big((c_{v}),(c^{\e\prime}_{v})\big)=\sum_{v\in
S}\text{inv}_{v}(c_{v}\lbe\cup c^{\e\prime}_{v}),
$$
where $\text{inv}_{\lbe v}\colon \text{Br}(F_{v})\ra\bq/\e\bz$ is
the usual invariant map of local class field theory. Now the dual of
\eqref{delta} is a map
$$
\delta^{D}\colon H^{\e 1}_{\rm{fl}}\be\big(U,\s
B^{\e\g^{\e\prime}}_{\ell^{r}}\big)\ra\displaystyle{\bigoplus_{v\in
S}}\,H^{\e 1}_{\rm{fl}}(F_{v},B_{ \ell^{r}}\!)
$$
and the following holds: if $c\in\prod_{\, v\in S} H^{\e
1}_{\rm{fl}}(F_{v},A_{ \ell^{r}}\!)$ and $x\in H^{\e
1}_{\rm{fl}}\be\big(U,\s B^{\e\g^{\e\prime}}_{\ell^{r}}\big)$, then
\begin{equation}\label{compatible2}
(c,\delta^{D}\be(x))=[\delta(c),x],
\end{equation}
where $[-,-]$ is the evaluation pairing \eqref{ev}.
Next, let
\begin{equation}\label{varrho}
\varrho\colon\displaystyle{\bigoplus_{v \in S}}\, H^{\e
0}_{\rm{fl}}(F_{v},A)\ra\displaystyle{\prod_{v\in S}} H^{\e
1}_{\rm{fl}}(F_{v},A_{\e \ell^{r}}\!)
\end{equation}
be the composite
$$
\displaystyle{\bigoplus_{v \in S}}\, H^{\e 0}_{\rm{fl}}(F_{v},A)\ra
\displaystyle{\bigoplus_{v\in S}} H^{\e 0}_{\rm{fl}}(F_{v},A)/\e
\ell^{r}\simeq \displaystyle{\prod_{v\in S}} H^{\e
0}_{\rm{fl}}(F_{v},A)/\e \ell^{r}\ra \displaystyle{\prod_{v\in S}}
H^{\e 1}_{\rm{fl}}(F_{v},A_{\e \ell^{r}}\!)
$$
(recall that $S$ is finite).

\begin{lemma} Let $c\in\prod_{\,v\in S}\, H^{\e
1}_{\rm{fl}}(F_{v},A_{\e \ell^{r}}\!)$. Then
$(c,\delta^{D}\be(x))=0$ for every
$$
x\in {\rm{Sel}}\big(U,\s B^{\e\g^{\e\prime}}\lbe\big)_{\ell^{
r}}\subset H^{\e 1}_{\rm{fl}}\be\big(U,\s
B^{\e\g^{\e\prime}}_{\ell^{r}}\big)
$$
if, and only if, $c=c_{1}+c_{2}$, with
$c_{1}\in\img\varrho$ and $c_{2}\in\krn\delta$.
\end{lemma}
\begin{proof}
The proof is similar to the proof of \cite{GA3}, lemma 5.5, using
\eqref{compatible2} and the commutative diagram
$$
\xymatrix{\displaystyle{\prod_{v \in S}} H^{\e 1}_{\rm{fl}}(F_{v},A_{ \ell^{r}})
\ar[rd]^{\delta} & & \\
\displaystyle{\bigoplus_{v \in S}}\, H^{\e 0}_{\rm{fl}}(F_{v},A)
\ar[u]^{\varrho}\ar[r] & H^{\e 1}_{\rm{fl}}(U,\s
B^{\e\g^{\e\prime}}_{ \ell^{r}})^{D} \ar@{->>}[r] &
\big({\rm{Sel}}\big(U,\s B^{\e\g^{\e\prime}}\lbe\big)_{\ell^{
r}}\big)^{D},}
$$
where $\delta$ and $\varrho$ are the maps \eqref{delta} and
\eqref{varrho}, respectively (the exactness of the bottom row of
this diagram follows from the definition of ${\rm{Sel}}\big(U,\s
B^{\g^{\prime}}\be\big)_{\ell^{ r}}$ and the local duality theorems
for abelian varieties).
\end{proof}
\begin{lemma} Assume that $\!\!\Sha^{1}(B)(\ell)$ is finite and let $a\in D^{\e 1}
\be\big(U,\s A^{\g}\big)$. If $a\in \ell^{\le r}H^{\e
1}_{\rm{fl}}\be\big(U,\s A^{\g}\lbe\big)$ and $\{a,a^{\e\prime}\}=0$
for every $a^{\e\prime}\in D^{\e 1}\be\big(U,\s
B^{\e\g^{\e\prime}}\lbe\big)_{\ell^{ r}}$, where $\{-,-\}$ is the
pairing \eqref{pair2}, then $a\in \ell^{r}D^{\e 1}\be\big(U,\s
A^{\g}\big)$.
\end{lemma}
\begin{proof} Consider the exact commutative diagram
$$
\xymatrix{ & H^{\e 1}_{\text{c}}\be\big(U_{\rm{fl}},\s A^{\g}\big) \ar[r]
\ar[d]^{ \ell^{r}} & H^{\e 1}_{\rm{fl}}\be\big(U,\s A^{\g}\big) \ar[d]^{ \ell^{r}} \\
\displaystyle{\bigoplus_{v\in S}}H^{\e 0}_{\rm{fl}}(F_{v},A)
\ar[r]^{\theta}\ar[d]^{ \varrho} & H^{\e
1}_{\text{c}}\be\big(U_{\rm{fl}},\s A^{\g}\big)
\ar[r]^{\psi}\ar[d]^{
\partial_{\e\text{c}}} & H^{\e 1}_{\rm{fl}}\be\big(U,\s A^{\g}\big) \ar[d]^{\partial}\\
\displaystyle{\prod_{v\in S}}H^{\e 1}_{\rm{fl}}(F_{v},A_{ \ell^{r}})
\ar[r]^{\delta^{\e\prime}} & H^{\e
2}_{\text{c}}\be\big(U_{\rm{fl}},\s A^{\g}_{\e \ell^{r}}\big) \ar[r]
& H^{\e 2}_{\rm{fl}}\be\big(U,\s A^{\g}_{\e \ell^{r}}\big)}
$$
where the middle row comes from \eqref{ccs} and the maps $\varrho$,
$\delta^{\e\prime}$ and $\partial_{\text{c}}$ are given by
\eqref{varrho}, \eqref{delta'} and  \eqref{partialc}, respectively.
Since $a\in D^{\e 1}\be\big(U,\s A^{\g}\e\big)=\img\psi$, there
exists $\widetilde{a}\in H^{\e 1}_{\text{c}}\be\big(U_{\rm{fl}},\s
A^{\g}\big)$ such that $\psi(\widetilde{a})=a$. Now, by hypothesis,
$0=\partial(a)=\partial(\psi(\widetilde{a}))$, whence
$\partial_{\text{c}}(\widetilde{a})=\delta^{\e\prime}(c)$ for some
$c\in\prod_{\,v\in S}H^{\e 1}_{\rm{fl}}(F_{v},A_{\e \ell^{r}}\!)$.
Now, if $x\in{\rm{Sel}}\big(U,\s B^{\e\g^{\e\prime}}\lbe\big)_{l^{
r}}$, then $\vartheta(x)\in D^{\e 1}\be\big(U,\s
B^{\e\g^{\e\prime}}\lbe\big)_{l^{ r}}$ by \eqref{selmer}, where
$\vartheta$ is the map \eqref{vartheta}. Therefore, by
\eqref{compatible2}, the definitions of $\delta$ and $\{-,-\}$ (see
\eqref{pair2} and \eqref{delta}) and \eqref{compatible},
\begin{eqnarray}
(c,\delta^{D}\be(x))\!\!\!&=&\!\!\![\e\delta(c),x\e]=
[\e\eta\e\delta^{\e\prime}(c),x]=[\e\eta\e
\partial_{\text{c}}(\widetilde{a}),x]\nonumber\\
\!\!\!&=&\!\!\!\langle\widetilde{a},\vartheta(x)\rangle=\{a,\vartheta(x)\}=0. \nonumber
\end{eqnarray}
Consequently, by Lemma 4.4, we may write $c\!=\!\varrho
(c^{\e\prime}_{\e 1})+ c_{\e 2}$ with $c^{\e\prime}_{\e
1}\in\prod_{\,v \in S}H^{\e 0}_{\rm{fl}}(F_{v},A)$ and $c_{\e
2}\in\krn\delta$. Now we have
\begin{eqnarray}
\eta\e\partial_{\e\text{c}}(\widetilde{a}-\theta(c^{\e\prime}_{\e 1}))\!\!\!&=&\!\!\!\eta\e\delta^{\e\prime}(c)-
\eta\e\delta^{\e\prime}\varrho(c^{\e\prime}_{\e 1})=\eta\e\delta^{\e\prime}(c-\varrho(c^{\e\prime}_{\e 1})) \nonumber \\
\!\!\!&=&\!\!\!\eta\delta^{\e\prime}(c_{\e 2})=\delta(c_{\e 2})=0.\nonumber
\end{eqnarray}
Thus the commutative diagram
$$
\xymatrix{H^{\e 1}_{\text{c}}\be\big(U_{\rm{fl}},\s A^{\g}\big)
\ar[r]^{\partial_{\text{c}}} \ar[d]& H^{\e
2}_{\text{c}}\be\big(U_{\rm{fl}},\s
A^{\g}_{\e \ell^{r}}\big) \ar[d]^{\eta} \\
\big(H^{\e 1}_{\rm{fl}}\be\big(U,\s B^{\e\g^{\e\prime}}\big)_{
\ell^{r}}\be\big)^{D}  \ar@{^{(}->}[r] & H^{\e
1}_{\rm{fl}}\be\big(U,\s B^{\e\g^{\e\prime}}_{\ell^{r}}\big)^{D} }
$$
shows that the element $\widetilde{a}-\theta(c^{\e\prime}_{\e 1})
\in H^{\e 1}_{\text{c}}\be\big(U_{\rm{fl}},\s A^{\g}\big)$ lies in
the kernel of the left-hand vertical map. Consequently, by
Proposition 4.3 (with the roles of $A$ and $B$ exchanged),
$\widetilde{a}-\theta(c^{\e\prime}_{\e 1})\in  \ell^{r}H^{\e
1}_{\text{c}}\be\big(U_{\rm{fl}},\s A^{\g}\big)$. Thus
$a=\psi(\widetilde{a})=\psi(\widetilde{a}-\theta(c^{\e\prime}_{\e
1}))\in \ell^{r}\img\psi= \ell^{r}D^{\e 1}\be\big(U_{\rm{fl}},\s
A^{\g}\big) $.
\end{proof}

\begin{lemma} If $\!\be\Sha^{1}(A)(\ell)$ is finite, then so also is $D^{1}\lbe\big(U
,\sa^{\g}\big)(\ell)$.
\end{lemma}
\begin{proof} Proposition 3.4 shows that the lemma is valid if $\varGamma=0$. Now the exact commutative diagram
\[
\xymatrix{\displaystyle\bigoplus_{v\in U } \g_{\be v}(k(v))\ar[r] & H^{1}_{\rm{fl}}(U,
\sao)\ar[d]
\ar[r] & H^{1}_{\rm{fl}}\big(U,\sa^{\g}\big)\ar[d]\ar[r] & \displaystyle
\bigoplus_{v\in U } H^{\e 1}_{\rm{fl}}(k(v),\g_{\lbe v})\\
&\displaystyle\prod_{v\in S}H^{\e 1}_{\rm{fl}}(F_{v},A)\ar@{=}[r] &
\displaystyle\prod_{v\in S} H^{\e 1}_{\rm{fl}}(F_{v},A),}
\]
whose top row is induced by \eqref{Gamma}, yields an exact sequence
$$
\displaystyle\bigoplus_{v\in U }\, \g_{\lbe v}(k(v))\ra
D^{1}(U,\sao) \ra
D^{1}\big(U,\sa^{\g}\e\big)\ra\displaystyle\bigoplus_{v\in U
}\,H^{\e 1}_{\rm{fl}}(k(v),\g_{\lbe v}).
$$
For any $i$, $H^{\e i}(k(v),\g_{\lbe v})$ is finite and equal to
zero for all but finitely many primes $v\in U $ \cite{S}, Chapter
XIII, \S1, Proposition 1, p.189. It follows that
$D^{1}\lbe\big(U,\sa^{\g}\big)(\ell)$ is finite if, and only if,
$D^{1}\lbe\big(U,\sao\big)(\ell)$ is finite. This completes the
proof.
\end{proof}

\begin{teorema} Let $\ell$ be any prime number such that $\!\!\Sha^{1}(A)(\ell)$
and $\!\!\Sha^{1}(B)(\ell)$ are finite. Then there exists a perfect
pairing of finite groups
$$
D^{1}\be\big(U,\s A^{\g}\big)(\ell)\times D^{1}\be\big(U,\s
B^{\e\g^{\e\prime}}\big)(\ell)\ra {\bq}/{\bz}.
$$
\end{teorema}
\begin{proof} The finiteness assertion is contained in the previous lemma.
Now there exists a canonical exact commutative diagram
$$
\xymatrix{0\ar[r]& D^{1}\lbe\big(U,\sa^{\g}\big)(\ell) \ar[r]
\ar[d]&
 H^{\e 1}_{\rm{fl}}\lbe\big(U,\sa^{\g}\big)(\ell)\ar[r] \ar[d] & \displaystyle
 \bigoplus_{v\in S}H^{\e 1}_{\rm{fl}}(F_{v},A)(\ell) \ar@{^{(}->}[d]\\
0\ar[r]& D^{\e 1}\be\big(U,\s B^{\e\g^{\e\prime}}\big)^{D} \ar[r] &
H^{\e 1}_{\text{c}}\be\big(U_{\rm{fl}},\s
B^{\e\g^{\e\prime}}\big)^{D}\ar[r] & \displaystyle\bigoplus_{v\in
S}H^{\e 0}_{\rm{fl}}(F_{v},B)^{D},}
$$
where the middle vertical map is induced by \eqref{pair'} and the
right-hand vertical map is injective by local duality. Proposition
4.3 implies that the kernel of the middle vertical map is contained
in $H^{\e 1}_{\rm{fl}}\lbe\big(U,\sa^{\g}\big)(\ell)_{\e
\ell-\text{div}}$. Now Lemma 4.5 and the finiteness of
$D^{1}\be\big(U,\s A^{\g}\big)(\ell)$ show that the left-hand
vertical map in the above diagram is injective. To complete the
proof, exchange the roles of $A$ and $B$.
\end{proof}
\begin{corollary} Assume that $\!\!\Sha^{1}(A)$ and $\!\!\Sha^{1}(B)$ are finite. Then
there exists a perfect pairing of finite groups
$$
D^{1}\be\big(U,\sao\big)\times D^{1}\be\big(U,\s
B\e\big)\ra{\bq}/{\bz}.
$$
\end{corollary}
\begin{proof} Take $\g=0$ and $\g^{\e\prime}=\Phi^{\e\prime}$ in the theorem and let
$\ell$ vary.
\end{proof}

\begin{teorema} Assume that $\!\!\Sha^{1}(A)$ and $\!\!\Sha^{1}(B)$ are finite. Then
there exists a perfect pairing of finite groups
$$
C_{A,\e F,\e S}\times C_{B,\e F,\e S}^{\e 1}\ra\bq/\bz,
$$
where $C_{B,\e F,\e S}^{\e 1}$ is the group \eqref{c1} associated to $B$.
\end{teorema}

\begin{proof} By the previous corollary and the existence of the perfect ``Cassels-Tate pairing" $\!\!\Sha^1(A)\times
\!\!\Sha^1(B)\ra\bq/\bz$ (\cite{HS}, Corollary 4.9, and
\cite{GA3}, Corollary 6.7), the dual of the second exact sequence appearing in Proposition 3.4 for $B$ is an exact sequence
$$
0\ra (C_{B,\e F,\e S}^{\e 1})^{D}\ra D^{1}\be\big(U,\sao\big)\ra
\Sha^1(A)\ra 0.
$$
The compatibility of the Cassels-Tate pairing with the pairings \eqref{pair} and \eqref{pair'} for $\g=0$ and $\g^{\e\prime}=\Phi^{\e\prime}$ (see \cite{GA-1}, Appendix) implies that the third map in the previous exact sequence is the same as that appearing in the first exact sequence of Proposition 3.4 for $A$. Thus there exists an exact commutative diagram
$$
\xymatrix{0\ar[r]& C_{A,\e F,\e S}\ar[r] &D^{1}\be\big(U,\sao\big)\ar[r] \ar@{=}[d] & \Sha^1(A)\ar@{=}[d]\ar[r]&0 \\
0\ar[r]&  (C_{B,\e F,\e S}^{\e 1})^{D}\ar[r] &
D^{1}\be\big(U,\sao\big)\ar[r] &
 \Sha^1(A)\ar[r]&0,}
$$
and this yields the result.
\end{proof}

\begin{remarks} (a) Assume that $\!\!\Sha^{1}(A)$, and therefore also $\!\!\Sha^{1}(B)$, is finite. Then Proposition 3.4 yields an exact sequence of finite groups
$$
0\ra C_{A,F,S}\ra D^{\e 1}(U,\sao) \ra D^{\e 1}(U,\s A\e)\ra
C_{A,F,S}^{\e 1} \ra 0
$$
whose dual, by Corollary 4.8 and the theorem, is the analogous exact sequence for $B$, i.e.,
$$
0\ra C_{B,F,S}\ra D^{\e 1}(U,\sbo) \ra D^{\e 1}(U,\s B\e)\ra
C_{B,F,S}^{\e 1} \ra 0.
$$

(b) Since $C_{A,\e F,\e S}$ is known to be finite, it is reasonable to question the necessity of the finiteness assumption in the theorem. In connection with this, W.McCallum \cite{McC}, Proposition 5.7, showed (independently of the finiteness assumption on $\!\!\Sha^{1}(A)$) the existence of a nondegenerate pairing
$$
H^{1}_{\rm{fl}}\be\big(U,\s A^{\g}\big)/{\rm{div}}\times
H^{1}_{\text{c}}\be\big(U_{\rm{fl}},\s
B^{\e\g^{\e\prime}}\big)_{\rm{tors}}\big/{\rm{div}}\ra {\bq}/{\bz}.
$$
See also \cite{Mi2}, Theorems III.3.7, p.317, and III.9.4, p.370.
Now, passing to the inverse limit over $r$ in the proof of Lemma 4.5 and using McCallum's result, one can show (again, independently of the finiteness assumption on $\!\!\Sha^{1}(A)$) that there exists a nondegenerate pairing of finite groups
$$
D^{1}\be\big(U,\s A^{\g}\big)(\ell)/\e{\text{$\ell$-div}}\times
D^{1}\be\big(U,\s
B^{\e\g^{\e\prime}}\big)(\ell)/\e{\text{$\ell$-div}}\ra {\bq}/{\bz}
$$
{\it provided} $\ell\neq p=\text{char}\e F$ in the function field
case (the problem is that the inverse limit over $r$ of the bottom
row of the big diagram in the proof of Lemma 4.5 might not be exact
if $\ell=p\,$\symbolfootnote[1]{Of course, we believe that the above
statement remains valid if $l=p$, but we see no obvious variant of
the proof of Lemma 4.5 that will prove this.}). However, even if the
above statement were verified for $\ell=p$ as well, one would still
need to check that $C_{A,\e F,\e S}\cap
D^{1}\be\big(U,\sao\big)_{\text{div}}=0$. So far, we have been
unable to make any progress on this problem.
\end{remarks}

\end{document}